\documentclass{article}
\newtheorem{df}{ \sc Definition}[section]

\newtheorem{ex}[df]{ \it Example}

\newtheorem{thrm}[df]{ \sc Theorem}

\newtheorem{re}[df]{ \it Remark}

\newtheorem{lem}[df]{ \sc Lemma}

\def\dim{{\rm dim}}

\def\ann{{\rm Ann}}
\def\socle{{\rm socle\ }}

\def\mpr#1{\;\smash{\mathop{\hbox to 20pt{\rightarrowfill}}\limits^{#1}}\;}
\def\epi#1{\;\smash{\mathop{\hbox to 20pt{\rightarrowfill}\hskip
-15pt\rightarrow}\limits^{#1\,}}\;}

\def\exseq#1#2#3{0\rightarrow #1 \rightarrow #2 \rightarrow #3 \rightarrow 0}
\def\O{{\mathcal O}}

\def\P{{\mathbb P}}

\def\F{{\mathcal F}}
\def\Hom{{\mathcal H om}}
\def\banica{{B\u anic\u a}}

\def\Q{{\mathbb Q}}

\def\Proof{\noindent\hskip4pt{\it Proof}.\ }
\def\qed{\hfill$\Box$\vskip10pt}

\input{amssymb.sty}
\usepackage{amscd}
\usepackage{amsmath}

\begin{document}

\author{Nicolae Manolache
}
\title{Linkage Extensions}
\date{}
\maketitle

\begin{abstract} Given two equidimensional Cohen-Macaulay local rings of the same
dimension, $A_1$ and $A_2$, we show that a simultaneous extension of $A_1$ by
the dualizing module of $A_2$ and of $A_2$ by the dualizing module of
$A_1$ is Gorenstein. This extends a theorem of Fossum (cf. \cite{Fo}).
The geometrical analogue of this is also considered. As an example, the pairs of double
lines in $\P ^3$ which are l.a.l linked are classified. This extends a
result of Migliore (cf. \cite{Mi}).

\end{abstract}

\section{Introduction}

In \cite{Fo}, generalizing a theorem of Reiten (cf. \cite{R}), Fossum
proves that the extension of a Cohen-Macaulay ring by a canonical module
is a Gorenstein ring. This idea is considered in a geometrical frame
by Ferrand (cf. \cite{Fe}), mainly for curves in a threefold. Both
these considerations can be interpreted in the frame of "algebraic
linkage" (cf. \cite{PS}).
In this paper we show that, in fact, if a ring extends two Cohen-Macaulay rings, each
by a dualizing module of the other, then this extension is Gorenstein.
We call such an extension "linkage extension". In a natural way one can
extend the notion to the case of schemes.

In general, given two (embedded) schemes, there is no (embedded) linkage
extension of them. In order to illustrate the "rarity" of linkage
extensions, we classify the pairs of double lines in $\P ^3$ which are
locally algebraically linked.

\section{Linked extension}

\subsection{The algebraic case}
All the rings considered here are commutative and noetherian. If $A$ is a local
ring, we denote by $m_A$ its maximal ideal.
\begin{df}{\label{d.e.}}Let $A_1$, $A_2$ be two rings, $M_1$ an $A_1$-module,
$M_2$ an $A_2$-module. We call  {\bf double extension of $A_1$ by $M_2$ and
$A_2$ by $M_1$} a couple of exact sequences:
$$
0\rightarrow M_2 \stackrel{i_2}{\rightarrow} B \stackrel{p_1}{\rightarrow} A_1
\rightarrow 0 \ ,
$$
$$
0\rightarrow M_1 \stackrel{i_1}{\rightarrow} B \stackrel{p_2}{\rightarrow} A_2
\rightarrow 0 \ ,
$$
where $B$ is a commutative ring, the maps $p_1$, $p_2$ are ring homomorphisms
and, for all
$b \in B$, $x_1 \in M_1$, $x_2 \in M_2$ we have:
\begin{eqnarray}
bi_1(x_1)=i_1(p_1(b)x_1) \nonumber \\
bi_2(x_2)=i_2(p_2(b)x_2) \nonumber
\end{eqnarray}
By abuse of language we call also {\bf double extension of $A_1$ by
$M_2$ and $A_2$ by $M_1$ } a ring $B$ which can be
inserted in a {\bf double extension of $A_1$ by $M_2$ and
$A_2$ by $M_1$} as above.
Also, when only $A_1$ and $A_2$ are given, we call {\bf double extension of
$A_1$ and $A_2$ } any ring $B$ with the property that it can be inserted in
extensions as above.
\qed
\end{df}

\begin{ex}
For each surjective map of rings $B\to B/a$ we realize canonicaly $B$ as
a double extension:
\begin{gather}
\exseq{a}{B}{B/a} \nonumber \\
\exseq{0:a}{B}{B/(0:a)} \nonumber
\end{gather}
\qed
\end{ex}
We shall consider a special case  of this notion:
\begin{df}{\label{l.e.}} When  $A_1$, $A_2$ in the above definition are
equidimensional Cohen-Macaulay rings of the same dimension, and $M_1$, $M_2$ are
respectively the dualizing modules $\omega_1$, $\omega _2$, a double extension
of $A_1$ by $\omega _2$ and of $A_2$ by $\omega _1$ is called
{\bf linkage extension of $A_1$ and $A_2$}. This means a couple of exact
sequences:
\begin{gather}
0\rightarrow \omega_2 \stackrel{i_2}{\rightarrow} B \stackrel{p_1}{\rightarrow}
A_1 \rightarrow 0 \ ,\nonumber \\
0\rightarrow \omega _1 \stackrel{i_1}{\rightarrow} B \stackrel{p_2}{\rightarrow}
A_2 \rightarrow 0 \ ,\nonumber
\end{gather}
where $B$ is a commutative ring, the maps $p_1$, $p_2$ are ring homomorphisms
and, for all
$b \in B$, $x_1 \in M_1$, $x_2 \in M_2$ we have:
\begin{eqnarray}
bi_1(x_1)=i_1(p_1(b)x_1) \nonumber \\
bi_2(x_2)=i_2(p_2(b)x_2) \nonumber
\end{eqnarray}
By abuse of language we call also a ring $B$ like above
a {\bf linkage extension of $A_1$ and $A_2$}
\qed
\end{df}

\begin{thrm}\label{gor}
Let $A_1$, $A_2$ be two equidimensional local Cohen-Macaulay local rings, of the same dimension.
In a linkage extension of $A_1$ and $A_2$ as above, the ring $B$ is Gorenstein.
\end{thrm}
\Proof The ring $B$ is Cohen-Macaulay of the dimension of $A_1$ and $A_2$. We shall
make an induction on this dimension. First two lemmas:
\begin{lem} In the situation from Definition \ref{l.e.} an element $b\in m_B$
is $B$-regular if and only if $p_1(b)$ is regular in $A_1$ and $p_2(b)$ is regular
in $A_2$.
\end{lem}
\Proof Everything follows from the observation that an element in a Cohen-Macaulay
ring is regular iff it is regular in the dualizing module.
\qed
\begin{lem} If $B$ is a linked extension of two artinian local rings $A_1$,
$A_2$, then
$$i_1({\socle} \omega _{A_1})=i_2({\socle} \omega _{A_2})={\socle}B\ .$$
\end{lem}
\Proof By symmetry, it is enough to show $i_1({\socle} \omega _{A_1})
={\socle}B$. For this, let $b \in {\socle }B$.
Then $bm_B=0$ and so $bi_1(\omega _{A_1})=0$. By definition this gives
$i_1(p_1(b)\omega _{A_1})=0$ and so $p_1(b)\omega _{A_1}=0$.
As $\ann_{A_1}\omega _{A_1}=0$ it follows $p_1(b)=0$. This proves the existence of an
element $b_2 \in \omega _{A_2}$ so that $b=i_2(b_2)$.
But $0=bm_B=i_2(b_2)m_B=i_2(b_2m_B)$ and so $b_2m_B=0$, i.e.
$b_2\in \socle \omega _{A_2}$. This shows $b\in i_2(\socle \omega _{A_2})$. Thus
 we proved
$\socle B \subset i_2(\socle \omega _{A_2})$. The other inclusion is evident.
\qed
We come back to the proof of the Theorem \ref{gor}. If $\dim B \ge 1$ let $b\in B$
be a regular element. Then $p_1(b)$ is regular in $A_1$ and $\omega _{A_1}$ and
$p_2(b)$ is regular in $A_2$ and $\omega _{A_2}$ and one shows easily that
$B/bB$ is a linked extension of $A_1/p_1(b)A_1$ and $A_2/p_2(b)A_2$. As $B$ is
Gorenstein iff $B/bB$ is Gorenstein, by repeating the above argument one reduces
the question to the artinian case. In this situation, as
$\omega _{A_i}\cong E(A_i/m_i)$, (the injective envelope of the residue field of $A_i$)
and the socle of $\omega_{A_i}$  is $A_i/m_i$, it follows $\socle B =A_i/m_i$ i.e.
the socle of $B$ is simple.
This proves that $B$ is Gorenstein. \qed
\begin{re}
When $A_1=A_2$, we obtain the Theorem of Fossum (cf. \cite{Fo}) that any
extension of a local Cohen-macaulay ring $A$ by its dualizing module  is Gorenstein.
We shall call such a ring {\bf a (Fossum) doubling of $A$}. In fact the proof
given here is "a splitting" of the original proof (cf. \cite{Fo}).
In a geometric frame the similar construction is known as {\bf Ferrand's
doubling}. This terminology is motivated by the fact that the multiplicity of
a doubling of a Cohen-Macaulay ring is double the multiplicity of that ring.
\qed
\end{re}
\begin{ex}
Let $R=k[x_1,\ldots ,x_n]$ with $n\ge 2$ and let $a_1$, $a_2$ be the ideals
$a_1=(x_1,x_2^2)$, $a_2=(x_1+x_2,x_2^2)$. Then $B:=R/(x_1^2+x_1x_2,x_2^2)$ is a
linkage extension of $A_1:=R/a_1$ and $A_2:=R/a_2$. The interest of this example
lies in the fact that $A_1$, $A_2$ are two "Fossum doublings" of $(A_1)_{\rm
red}=(A_2)_{\rm red}=k[x_3,\ldots x_n]=:A$ and $B$ is a multiplicity $4$
structure on $A$ which is not a doubling of $A_1$ or $A_2$. In fact $B$ can be
realized as a doubling of $A_3=R/(x_1^2,y)$. It can be realized also as a
linkage extension of $A$ and a triple structure on $A$.
\qed
\end{ex}

\subsection{The Geometric Case}
\begin{df}\label{geom-ext}
If $X_1$, $X_2$ are two algebraic schemes, we call {\bf double extension}
of them a couple of exact sequences:
$$
0\rightarrow \F _2 \stackrel{i_2}{\rightarrow} \O _Y
\stackrel{p_1}{\rightarrow} \O _{X_1}
\rightarrow 0 \ ,
$$
$$
0\rightarrow  \F _1 \stackrel{i_1}{\rightarrow} \O _Y
\stackrel{p_2} {\rightarrow} \O_{X_2}
\rightarrow 0 \ ,
$$
where $\F_i$ are $\O_{X_i}$-modules and $Y$ is an algebraic scheme,
such that locally one has a double extension in the sense of
\ref{d.e.}

If $X_1$, $X_2$ are equidimensional locally Cohen-Macaulay schemes,
of the same dimension, we call {\bf linkage extension} of $X_1$ and $X_2$
a double extension of the shape:

\begin{gather}
0\rightarrow \omega_2 \otimes L_2 \stackrel{i_2}{\rightarrow} \O _Y
\stackrel{p_1}{\rightarrow} \O _{X_1}
\rightarrow 0 \ , \label{ex1}\\
0\rightarrow \omega _1 \otimes L_1 \stackrel{i_1}{\rightarrow} \O _Y
\stackrel{p_2} {\rightarrow} \O_{X_2}
\rightarrow 0 \ ,
\end{gather}
where , $L_i$ are invertible $\O_{X_i}$-modules and
$\omega _i$ are the corresponding dualizing sheaves.
By abuse of language we say also that $Y$ is a {\bf linkage
extension} of $X_1$, $X_2$.
\qed
\end{df}

\begin{thrm}\label{gor-sch}
If $X_1$, $X_2$ are equidimensional locally Cohen-Macaulay schemes,
of the same dimension, any linkage extension of them is locally
Gorenstein.
\end{thrm}
\Proof The question is local, so \ref{gor} applies.
\qed

\begin{lem}
If $X$ is a linkage extension of two equidimensional locally Cohen-Macaulay
schemes $X_1$, $X_2$ of the same dimension, then, with the notation from
the definition \ref{geom-ext}:
$$
L_1=\omega _Y^{-1} |X_1, \ \ L_2=\omega _Y^{-1} |X_2 \ .
$$
\end{lem}
\Proof

Applying the dualizing functor $\Hom (?, \omega _Y)$ to the exact
sequence (\ref{ex1}) and then tensoring with $\omega _Y^{-1}$ one gets:
$$
\exseq{\omega_1\otimes \omega _Y^{-1}}{\O_Y}{\O_{X_2}\otimes
L_2^{-1}\otimes \omega_Y^{-1}}
$$
Here $O_{X_2}\otimes L_2^{-1}\otimes \omega_Y^{-1}$ is the
structural sheaf of a closed subscheme of $Y$ which
coincides locally with $X_2$, so it is $\O _{X_2}$. Then
$\omega_1\otimes \omega _Y^{-1}=\omega_1\otimes L_1$. This shows
$L_1\cong \omega _Y^{-1} |X_1$. Analogously $L_2\cong \omega _Y^{-1} |X_2$.
\qed

\begin{re}
If, given  $X_1$, $X_2$, there exists a $Y$ which is  a linkage
extension of  $X_1$, $X_2$, we say also that $X_1$ and $X_2$ are
locally Gorenstein linked. When $Y$ is locally complete intersection
we say that $X_1$ and $X_2$ are {\bf locally algebraically linked,
l.a.l. for short},
(cf. \cite{M2}, where this terminology was introduced, inspired by
the notion of {\bf algebraic linkage} of \cite{PS}).
\qed
\end{re}

If $X_1$ and $X_2$ are both embedded in a (let say smooth) scheme
$P$ and we require that the linkage extension $Y$ to be also
closed subscheme of $P$, we say that {\bf $Y$ is an embedded (in
$P$) linkage extension of $X_1$ and $X_2$}. In general, given $X_1$
and $X_2$ in $P$ there is no linkage extension (in $P$) of them.

The aim of the next theorem is to classify the pairs of double lines
in $\P^3$ which are locally algebraically  linked.
\begin{thrm}
If two double lines $Y_1$, $Y_2$ in $\P ^3$ with supports respectively $X_1$,
$X_2$,  are l.a.l. then they are in one of
the  following situation:

\begin{description}
\item (i) $Y_1$, $Y_2$ are disjoint

\item (ii) in convenient homogeneous coordinates $(x:y:z:u)$ (i.e.
after an automorphism of $\P^3$) they are defined by ideals of the form:
\begin{gather}
I_{Y_1}=(ax+by,x^2,xy,y^2) \nonumber \\
I_{Y_2}=(cx+dz,x^2,xz,z^2) \ ,\nonumber
\end{gather}
where $a(z,u)$, $b(z,u)$ are homogeneous forms in $z,u$ of the same degree
$r_1$, without common zeros on $X_1$, $c(y,u)$, $d(y,u)$ are homogeneous forms
in $y,u$ of the same degree $r_2$, without common zeros on $X_2$,
such that:
\begin{description}
\item (a) $b(0:1)\neq 0,\ d(0:1)\neq 0$,
\item or
\item (b) $b(0:1)=0 ,\ d(0:1)=0$ and $a(0:1)\displaystyle\frac{\partial b}{\partial
 z}(0:1)= c(0:1)\frac{\partial d}{\partial y}(0:1)$
 \end{description}

\item (iii) if $X_1=X_2$ then, either

 \begin{description}
 \item  (a) $Y_1=Y_2$,
 \item or
 \item  (b) in convenient homogeneous coordinates $(x:y:z:u)$ they are defined
by ideals of the form:
\begin{gather}
I_{Y_1}=(ax+by,x^2,xy,y^2) \nonumber \\
I_{Y_2}=(ax-by,x^2,xy,y^2) \ ,\nonumber
\end{gather}
where $a(z,u)$, $b(z,u)$ are homogeneous forms in $z,u$ of the same degree
$r$.
\end{description}
\end{description}
\end{thrm}
\begin{Proof}
We have to show that, when $X_1$ and $X_2$ meet in a point,
(ii) is fulfilled and that, when $X_1=X_2$,
(iii) is  fulfilled.

If $X_1$, $X_2$ have a point in common, one may suppose that the
equations of $X_1$ and $X_2$ are $x=y=0$, respectively
$x=z=0$ in convenient homogeneous coordinates $(x:y:z:u)$. Then
double structures $Y_1$, $Y_2$ on them are given by ideals like those in (ii).
We have to determine when $Y_1\cup Y_2$ is a complete intersection
in $X_1\cap X_2=(0:0:0:1)$. In the affine space $u\ne 0$, with
coordinates $\xi=\displaystyle\frac{x}{u}$,
$\eta=\displaystyle\frac{y}{u}$, $\zeta=\displaystyle\frac{z}{u}$,
one has :
\begin{eqnarray}
I_{Y_1}=(\alpha \xi +\beta \eta ,\ \xi ^2 ,\ \xi\eta ,\ \eta
^2) \nonumber \\
I_{Y_2}=(\gamma \xi +\delta \zeta ,\ \xi ^2 ,\ \xi\zeta ,\ \zeta
^2) \ ,\nonumber
\end{eqnarray}
where $\alpha $, $\beta $ are functions in $\zeta$ and $\gamma$,
$\delta $ are functions in $\eta $.

 A direct computation shows that,
 if only one of the $\beta (0)$, $\delta (0)$ is $0$, then $Y_1\cup
 Y_2$ is not a complete intersection in $Y_1\cap Y_2$.

If $\beta (0)\ne 0$, $\delta (0)\ne 0$ and with new parameters $Y=
 \alpha \xi +\beta \eta$ and $Z=\gamma \xi +\delta \zeta$, the
 ideals are locally: $I_{Y_1}=(Y,\ \xi ^2)$, $I_{Y_2}=(Z,\ \xi ^2)$ and
 $I_{Y_1\cup Y_2}=(YZ,\ \xi ^2)$. This hives the case (a).

 Consider $\beta (0)=0$, $\delta (0)=0$, i.e. $\beta =\beta _1\zeta$,
 $\delta =\delta _1 \eta$. Then $\alpha (0)\ne 0$ and
 $\gamma (0)\ne 0$. The ideal of $Y_1\cup Y_2$ is then
 the intersection :
 \begin{eqnarray}
 (\alpha (0)\xi +\beta _1 (0)\zeta\eta+L\xi\zeta +M\zeta ^2\eta,\
 \xi^2,\ \xi\eta,\ \eta^2) \ \cap \nonumber  \\
 (\gamma (0)\xi +\delta _1(0)\eta\zeta+N\xi\eta +P\eta ^2\zeta,\
 \xi^2,\ \xi\zeta,\ \zeta^2) \ ,\nonumber
 \end{eqnarray} i.e. it is
 $$((\alpha (0)\xi +\beta _1 (0)\zeta\eta,\
 \xi^2,\ \xi\eta,\ \eta^2)\cap (\gamma (0)\xi +\delta _1
 (0)\eta\zeta+N\xi\eta +P\eta ^2\zeta,\
 \xi^2,\ \xi\zeta,\ \zeta^2) \ .$$
 This ideal defines a complete intersection iff there exists a
 $\lambda \ne 0$ such that $\alpha (0)\xi +\beta _1 (0)\zeta\eta=\lambda
 (\gamma (0)\xi +\delta _1(0)\eta\zeta )$. This condition translates
 to (b).

 Consider now the case when $X_1=X_2=:X$. If $Y_1$ and
 $Y_2$ are l.a.l., then there exists a multiplicity $4$ structure $Y$  on
 $X$ which is a linked extension of $Y_1$ and $Y_2$.
 Consider first the case
 $Y$ is a quasiprimitive structure in the sense of {\banica}  and Forster
 (cf. \cite{BF2} or \cite{M3}).
 Then $Y$  is a doubling of a doubling $Y'$ of $X$ (cf. \cite{M3}, {\bf Lemma
 2.10}). The curve $Y'$ is obtained canonically from $Y$, being a member
 of the filtrations introduced in \cite{BF2}, \cite{M2} and \cite{M3},
 which coincide for quasiprimitive structures. In the \banica - Forster
 filtration $Y'$ is obtained throwing away the embedded points of $Y\cup
 X^{(2)}$. As $Y_1$ and $Y_2$ are both contained in $X^{(2)}$, they should
 be contained in $Y'$, i.e. one has $Y_1=Y_2=Y'$.

 If $Y$ is not a quasiprimitive structure, then according to \cite{BF1} or
 \cite{BF2} it is globally complete intersection, and, in convenient
 coordinates, its ideal is $I_Y=(x^2, y^2)$. In these coordinates, with
 convenient forms $a(z,u)$, $b(z,u)$ of the same degree, without common
 zeros along $X$, the ideal of $X_1$ is $I_1=(ax+by,x^2,xy,y^2)$.
 Then the ideal of $Y_2$ should be $I_Y:I_1 = (ax-by,x^2,xy,y^2)$.
 \qed
 \end{Proof}

{\bf Acknowledgement.} The author was partially supported by {\em
Contract\\
CERES 3-28, 2003-2005 and by Contract CERES 2004-2006}. Thanks are due also
to the Institute of Mathematics of the University of Oldenburg, especially
to Prof. U. Vetter, for the warm hospitality in decmber 2004, when last
touches to this paper were done.

\noindent Nicolae Manolache\\
Institute of Mathematics "Simion Stoilow"\\
of the Romanian Academy \\
P.O.Box 1-764
Bucharest, RO-014700

\noindent e-mail: nicolae.manolache@imar.ro

\begin{thebibliography}{HVDV}
\small
\def\by{ \sc}
\def\book{~: \rm}
\def\paper{~: \em }
\def\jour{~, \rm}
\def\inbook{~, \rm}
\def\publ{\rm, }
\def\vol{~, \bf}
\def\no{~, \rm no.~}
\def\yr#1{ \rm (#1)}
\def\pages{~, \rm p.~}
\def\endref{\rm}

\bibitem[BF1]{BF1}\by C.~B\u anic\u a, O.~Forster
\paper Sur les Structures Multiples (manuscript)
\yr{1981}\endref

\bibitem[BF2]{BF2}\by C.~B\u anic\u a, O.~Forster
\paper Multiplicity structures on Space Curves\jour Contemporary
mathematics \vol 58\yr {1986}\endref

\bibitem[BM]{BM}\by C.~B\u anic\u a, N.~Manolache
\paper Rank 2 Stable vector bundles on ${\Bbb P}^3(\Bbb C)$\jour
Math. Z. \vol 190\yr
{1985}\pages 315-339\endref

\bibitem[BM']{BM'} \by C.~B\u{a}nic\u{a}, N.~Manolache \paper
Remarks on rank 2 stable vector bundles on ${\Bbb P}^3$ with Chern classes
$c_1=-1$, $c_2=4$\jour Preprint INCREST\vol 104\yr {1981}\endref

\bibitem[BM"]{BM"} \by C.~B\u{a}nic\u{a}, N.~Manolache  \paper
Moduli space $M(-1,4)$: Minimal spectrum\jour Preprint
INCREST\vol19\yr {1983} \endref

\bibitem[Ba]{Ba}\by H.~Bass \paper On the ubiquity of Gorenstein rings\jour
Math. Z.\vol 82\yr{1963}\pages 8-28 \endref

\bibitem[BE]{BE}\by D.~Bayer, D.~Eisenbud \paper Ribbons and their
canonical  embeddings
\jour Trans. Amer. Math. Soc.\vol 347\yr{1995}\pages 757-765\endref

\bibitem[Fe]{Fe}\by D.~Ferrand
\paper Courbes Gauches et Fibr\' es de Rang 2\jour C.R. Acad. Sci. Paris\vol 281\yr
{1975}\pages 345-347\endref

\bibitem[Fo]{Fo}\by R.~Fossum \paper Commutative Extensions by Canonical
Modules are Gorenstein Rings \jour Proc. of the Amer. Math. soc.\vol
40\yr{1973}\pages 395-400\endref

\bibitem[H1]{H1}\by R.~Hartshorne
\paper Stable Vector Bundles of Rank 2 on $\P^3$ \jour Math. Ann.\vol 238\yr
{1978}\pages 229-280\endref

\bibitem[H2]{H2}\by R.~Hartshorne
\paper Stable reflexive Sheaves of Rank 2 on $\P^3$\jour Math. Ann.\vol 254\yr
{1978}\pages 121-176\endref

\bibitem[M1]{M1}\by N.~Manolache
\paper Rank 2 Stable Vector Bundles on $\P^3$with Chern classes $c_1=-1$, $c_2=2$
\jour Rev. Roumaine Math. pures et Appl.\vol 26\yr{1981}\pages 1203-1209\endref

\bibitem[M2]{M2}\by N.~Manolache
\paper Cohen-Macaulay Nilpotent Structures
\jour Rev. Roumaine Math. pures et Appl.\vol 31\yr{1986}\pages 563-575\endref

\bibitem[M3]{M3}\by N.~Manolache
\paper Codimension Two Linear Varieties with Nilpotent Structures
\jour Math. Z.\vol 210\yr{1992}\pages 573-580\endref

\bibitem[M4]{M4}\by N.~Manolache
\paper Curves $Y$ of Degree 6 with $\omega _Y=\O _Y(-1)$ on a Smooth Quadric $\Q_3$
\jour Annali di Matematica pura ed  applicata (IV)\vol 167\yr
{1994}\pages 226-241\endref

\bibitem[M5]{M5} \by N. Manolache \paper Multiple Structures on Smooth Support
\jour Math. Nachr.\vol 167\yr{1994}\pages 157-202 \endref

\bibitem[M6]{M6} \by N.~Manolache \paper Double Rational Normal Curves
with Linear Syzygies\jour
Manuscripta Mathematica \vol 104, No. 4\yr{2001}\pages 503-517

\bibitem[Mi]{Mi} \by J. Migliore \paper On linking double lines \jour
Trans. Am. Math. soc., \vol 294 \yr{1986} \pages 177-185

\bibitem[NNS]{NNS} \by U.~Nagel, R.~Notari, M.~L.~Spreafico\paper The Hilbert Scheme of
Degree Two Curves and Certain Ropes\jour arXiv:math.AG/0311537

\bibitem[N]{N}\by S.~Nollet \paper The Hilbert Schemes of Degree Three Curves
\jour Ann. Scient. E.N.S., $4^e$ S\' e rie\vol 30\yr{1997}\pages 367-384

\bibitem[NS]{NS}\by S.~Nollet, E. Schlesinger\paper Hilbert Schemes of
Degree Four Curves\jour arXiv:math.AG/0112167

\bibitem[OSz]{OSz}\by G.~Ottaviani, M.~Szurek \paper On Moduli of Stable
2-Bundles with Small Chern Classes on $Q_3$ \jour Annali di Matematica pura ed
applicata (IV),\vol 167\yr{1994}\pages 191-241\endref

\bibitem[PS]{PS}\by C.~Peskine, L.~Szpiro \paper Liaison des Vari\' etes Alg\' ebriques
\jour Invent. Math.\vol 26\yr{1977}\pages 271-302\endref

\bibitem[R]{R}\by I.~Reiten \paper The Converse to a Theorem of Sharp on
Gorenstein Modules\jour Proc. Amer. math. Soc.\vol 32\yr{1972}\pages
417-420\endref

\bibitem[Sz]{Sz}\by L.~Szpiro \paper Equations Defining Space Curves
\jour Tata institute Lect. Notes (notes by N. Mohan Kumar, Bombay\yr{1979} \endref

\bibitem[V1]{V1} \by J.~E.~Vatne \paper Towards a Classification of Multiple Structures
\jour Ph.D. thesis,  University of Bergen\yr {2001} \endref

\bibitem[V2]{V2} \by J.~E.~Vatne \paper Multiple Structures\jour
arXiv:math.AG/0210042 \endref

\bibitem[V3]{V3} \by J.~E.~Vatne \paper Monomial Multiple Structures\jour
arXiv:math.AG/0210101 \endref

\end{thebibliography}
\end{document}